\documentclass{amsart}

\usepackage{amssymb,bm,amsthm,color,scalerel}
\usepackage{hyperref}

\usepackage{tikz}
\usetikzlibrary{positioning}
\usetikzlibrary{svg.path}
\definecolor{orcid_color}{HTML}{A6CE39}

\hypersetup{pdfborder={0 0 0}} % Remove the border of a hyperlink
\DeclareRobustCommand{\orcidicon}{%
	\raisebox{.2mm}{\scalerel*{%
	\begin{tikzpicture}[xscale=1,yscale=-1,transform shape]
	\filldraw[color=orcid_color] svg {M256,128c0,70.7-57.3,128-128,128C57.3,256,0,198.7,0,128C0,57.3,57.3,0,128,0C198.7,0,256,57.3,256,128z};
	\filldraw[color=white] svg {M86.3,186.2H70.9V79.1h15.4v48.4V186.2z} svg {M108.9,79.1h41.6c39.6,0,57,28.3,57,53.6c0,27.5-21.5,53.6-56.8,53.6h-41.8V79.1z M124.3,172.4h24.5
		c34.9,0,42.9-26.5,42.9-39.7c0-21.5-13.7-39.7-43.7-39.7h-23.7V172.4z} svg {M88.7,56.8c0,5.5-4.5,10.1-10.1,10.1c-5.6,0-10.1-4.6-10.1-10.1c0-5.6,4.5-10.1,10.1-10.1
		C84.2,46.7,88.7,51.3,88.7,56.8z};
	\end{tikzpicture}}{|}}%
}
\newcommand{\orcid}[1]{\href{https://orcid.org/#1}{\orcidicon}}

\theoremstyle{plain}
\newtheorem{theorem}{Theorem}[section]
\theoremstyle{definition}
\newtheorem{example}[theorem]{Example}
\newtheorem*{remark}{Remark}

\newcommand{\arxiv}[1]{arXiv:\href{http://arxiv.org/abs/#1}{#1}}

\begin{document}

\title{The real equiangular tight frames obtained from rank $3$ graphs}
\author[Ei. Bannai]{Eiichi Bannai}
\address{Faculty of Mathematics, Kyushu University (emeritus), Japan}
\email{bannai@math.kyushu-u.ac.jp}
\author[Et. Bannai]{Etsuko Bannai}
\email{et-ban@rc4.so-net.ne.jp}
\author[C.-Y. Lee]{Chin-Yen Lee}
\address{Department of Mathematics, National Central University, Taoyuan, 32056, Taiwan}
\email{109281001@cc.ncu.edu.tw}
\author[H. Tanaka]{Hajime Tanaka\,\orcid{0000-0002-5958-0375}}
\address{Research Center for Pure and Applied Mathematics, Graduate School of Information Sciences, Tohoku University, Sendai 980-8579, Japan}
\email{htanaka@tohoku.ac.jp}
\author[W.-H. Yu]{Wei-Hsuan Yu\,\orcid{0000-0003-3569-7795}}
\address{Department of Mathematics, National Central University, Taoyuan, 32056, Taiwan}
\email{whyu@math.ncu.edu.tw}

 % 42C15: General harmonic expansions, frames; 20B15: Primitive groups; 05E30: Association schemes, strongly regular graphs; 05C50: Graphs and linear algebra (matrices, eigenvalues, etc.)
\begin{abstract}
We present all nontrivial real equiangular tight frames $\{\varphi_m\}_{m=1}^M$ in $\mathbb{R}^N$ obtained as spherical embeddings of primitive rank $3$ graphs on $M$ vertices, and those such that one of their associated $M$ strongly regular graphs on $M-1$ vertices is a primitive rank $3$ graph.
\end{abstract}

\maketitle

\hypersetup{pdfborder={0 0 1}} % Put back the borders of hyperlinks

%%%%%%%%%%%%%%%%%%%%%
%%%%%%%%%%%%%%%%%%%%%
\section{Introduction}

It is well known that a set $\{ \varphi_m \}_{ m = 1 }^M$ of $M\,(\geqslant N)$ unit vectors in $\mathbb{R}^N$ satisfies
\begin{equation*}
	\max_{ m \ne n } | \langle \varphi_m ,\varphi_n \rangle | \geqslant \sqrt{ \frac{ M - N }{ N ( M - 1 ) } },
\end{equation*} 
and equality is attained if and only if $\{ \varphi_m \}_{ m = 1 }^M$ is \emph{equiangular}, i.e., $| \langle \varphi_m ,\varphi_n \rangle |$ is constant for $m\ne n$, and it is a \emph{tight frame}, i.e., there is a constant $c > 0$ such that
\begin{equation*}
	\sum_{ m = 1 }^M \langle \varphi, \varphi_m \rangle^2 = c \| \varphi \|^2
\end{equation*}
for all $\varphi \in \mathbb{R}^N$; see, e.g., \cite{Welch1974IEEE}, \cite[\S10.6.2]{BH2012B}.
It seems that the concept of real equiangular tight frames (real ETFs) first appeared in Van Lint and Seidel \cite{LS1966IM}; see also Lemmens and Seidel \cite{LS1973JA}.
We refer the reader to \cite{CRT2008P,FJMP2018JCTA,FM2015pre,FMT2012LAA,IJM2020DCG,IM2019pre,Singh2010LAA,Waldron2009LAA} for some of the recent activities on ETFs.
This paper presents all nontrivial real ETFs obtained as spherical embeddings of primitive rank $3$ graphs on $M$ vertices, and those such that one of their associated $M$ strongly regular graphs on $M-1$ vertices is a primitive rank $3$ graph.

This research was started when Ei.B. received an email from Bob Griess dated July 29, 2012.
The content of his email was essentially as follows.

\begin{quote}
``There are interesting examples of spherical codes which come from VOA\footnote{VOA stands for `vertex operator algebra.'} theory.
For instance, there is a code in $156$-dimensional Euclidean space with $496$ unit vectors and group $O^+_{ 10 } ( 2 )$.
If $x$ is one unit vector, it has inner product $0$ with $255$ others and inner product $1/8$ with $240$ others.
(The points of the spherical code are rescales of the so-called conformal vectors of central charge $1/2$).
Do you know anything about this code?
For example, is it close to being extremal in any sense?''
\end{quote}

\noindent
(He also mentioned that the above procedure is applied to the Monster simple group acting in $196884$ dimensions.)
The reply to his inquiry was essentially as follows. 

\begin{quote}
``The $496$ points are on the hyperplane in $\mathbb{R}^{156}$ which is at distance $1/4$ from the origin.
If we regard this $496$ point set as a spherical code on the unit sphere in the $155$-dimentional Euclidean space, then the inner products become $1/15$ and $-1/15$.
So, we have an equiangular line system.
Moreover, this gives an extremal property of being a real equiangular tight frame. 
Also, this construction of real equiangular tight frames is generalized for $O^+_{ 2 m } ( 2 )$ and $O^-_{ 2 m } ( 2 )$ (for all $m$) acting as transitive rank $3$ groups (i.e., strongly regular graphs) on the set of non-isotropic points.''
\end{quote}

Then, Ei.B., Et.B., and H.T. quickly found four other infinite families by looking at the eigenmatrices of various primitive rank $3$ graphs but did not publish the results.
In Winter 2021, C.-Y.L. and W.-H.Y. joined the project of finding \emph{all} nontrivial real ETFs obtained as spherical embeddings of primitive rank $3$ graphs.
The classification of primitive rank $3$ graphs was the culmination of efforts by many researchers such as Foulser, Kallaher, Kantor, Liebler,  Liebeck, and Saxl.
Their results were inevitably scattered in a mass of papers and were not easy to follow.
Our work was then greatly facilitated by the recent book written by Brouwer and Van Maldeghem \cite{BVM2022B}, which includes descriptions of all primitive rank $3$ graphs.
Besides the above $2+4=6$ infinite families, it turns out that there are only two sporadic examples.
Associated with every nontrivial real ETF $\{ \varphi_m \}_{ m = 1 }^M$ are $M$ strongly regular graphs on $M-1$ vertices having the same parameters.
As a related problem, we also determine all nontrivial real ETFs such that one of these strongly regular graphs is a primitive rank $3$ graph.
Our main results are Theorems \ref{classification} and \ref{rank 3 Waldron}.

This paper is organized as follows.
Section \ref{sec: SRGs} recalls the necessary facts about strongly regular graphs and their spherical embeddings.
In Section \ref{sec: rank 3 graphs}, we summarize the classification of primitive rank $3$ graphs based on Brouwer and Van Maldeghem \cite{BVM2022B} and identify all those whose spherical embeddings give rise to real ETFs.
Section \ref{sec: real ETFs from SRGs} describes these six infinite families and two sporadic graphs.
Finally, Section \ref{sec: Waldron ETFs} discusses nontrivial real ETFs, one of whose associated strongly regular graphs is a primitive rank $3$ graph.

%%%%%%%%%%%%%%%%%%%%%
%%%%%%%%%%%%%%%%%%%%%
\section{Spherical embeddings of strongly regular graphs}
\label{sec: SRGs}

Let $\Gamma$ be a finite simple graph with vertex set $X$, that is neither complete nor edgeless.
Recall that $\Gamma$ is called a \emph{strongly regular graph} (\emph{SRG}) \emph{with parameters} $(v,k,\lambda,\mu)$ if it has $v$ vertices and is $k$-regular, and every pair of adjacent (resp.~distinct and nonadjacent) vertices has exactly $\lambda$ (resp.~$\mu$) common neighbors.
For the rest of this section, suppose that $\Gamma$ is an SRG with parameters $(v,k,\lambda,\mu)$.
The complement $\overline{\Gamma}$ of $\Gamma$ is an SRG with parameters $(v,\bar{k},\bar{\lambda},\bar{\mu})$, where $\bar{k}=v-k-1$, $\bar{\lambda}=v-2k+\mu-2$, and $\bar{\mu}=v-2k+\lambda$; cf.~\cite[\S 1.1.2]{BVM2022B}.
We also assume that $\Gamma$ is \emph{primitive}, i.e., $\Gamma$ is neither a complete multipartite graph $(\mu=k)$ nor a union of complete graphs $(\mu=0)$.
Let $A$ and $\overline{A}=J-I-A$ be the adjacency matrices of $\Gamma$ and $\overline{\Gamma}$, respectively, where $J$ is the all one's matrix.
Then $\bm{A} := \mathrm{span} \{ I, A, \overline{A} \}$ is a three-dimensional algebra since $A^2=kI+\lambda A+\mu \overline{A}$.
The algebra $\bm{A}$ is known as the \emph{Bose--Mesner algebra} of $\Gamma$; cf.~\cite[\S1.3.2]{BVM2022B}.
The matrix $A$ has eigenvalue $k$ with multiplicity one since $\Gamma$ is connected and $k$-regular.
It has two other distinct eigenvalues $r>s$ and these are the roots of the quadratic equation $\xi^2+(\mu-\lambda)\xi+(\mu-k)=0$; cf.~\cite[\S 1.1.1]{BVM2022B}.
Let $E_k, E_r, E_s$ be the orthogonal projections onto the eigenspaces $V_k, V_r, V_s$ of $A$ with eigenvalues $k, r, s$, respectively.
We note that $E_k = v^{ - 1 } J$ and that $E_k, E_r$, and $E_s$ form a basis of $\bm{A}$.
Define the $3\times 3$ matrices $P$ and $Q$ by $( I, A, \overline{A} ) = ( E_k, E_r, E_s ) P$ and $v ( E_k, E_r, E_s ) = ( I, A, \overline{A} ) Q$.
We call $P$ and $Q$ the \emph{first} and \emph{second eigenmatrices} of $\Gamma$, respectively.
We have
\begin{equation*}
	P = \begin{pmatrix} 1 & k & \bar{k} \\ 1 & r &\bar{s} \\ 1 & s & \bar{r} \end{pmatrix}\!,
\end{equation*}
where $\bar{r}=-s-1$ and $\bar{s}=-r-1$ are the distinct eigenvalues of $\overline{A}$ other than $\bar{k}$.
Moreover, it follows that
\begin{equation*}
	Q = \begin{pmatrix} 1 & f & g \\ 1 & fr/k & gs/k \\ 1 & f\bar{s}/\bar{k} & g\bar{r}/\bar{k} \end{pmatrix}\!,
\end{equation*}
where $f=\dim V_r$ and $g = \dim V_s=v-f-1$.
Observe in particular that
\begin{equation}\label{E_r}
	E_s = \frac{ g }{ v } \left( I + \frac{ s }{ k } A + \frac{ \bar{r} }{ \bar{k} } \overline{A} \right).
\end{equation}

For every $x \in X$, let $( E_s )_x$ denote the $x^{ \mathrm{th} }$ column of $E_s$, and consider the set $\{ \varphi_x : x \in X \}$ of $v$ (column) vectors in the Euclidean space $V_s \, ( \cong \mathbb{R}^g )$ by
\begin{equation*}
	\varphi_x = \sqrt{ \frac{ v }{ g } } ( E_s )_x \qquad ( x \in X ).
\end{equation*}
Then it follows from \eqref{E_r} that the $\varphi_x$ are unit vectors in $V_s$ and that
\begin{equation}\label{angles}
	\langle \varphi_x, \varphi_y \rangle = \begin{cases} s/k & \text{if $x, y$ are adjacent}, \\ \bar{r}/\bar{k} & \text{if $x, y$ are nonadjacent}, \end{cases} \qquad ( x, y \in X, \ x \ne y ).
\end{equation}
In this paper, we call $\{ \varphi_x : x \in X \}$ the \emph{spherical embedding} of $\Gamma$.
Since
\begin{equation*}
	\sum_{x\in X} \langle \varphi, \varphi_x \rangle^2 = \frac{v}{g} \varphi^{\mathsf{T}} E_s \varphi  = \frac{v}{g} \| \varphi \|^2 \qquad (\varphi\in V_s),
\end{equation*}
the spherical embedding of $\Gamma$ gives a real ETF with $(M,N)=(v,g)$ if and only if it is equiangular, i.e.,
\begin{equation}\label{equiangular}
	\frac{s}{k}=-\frac{\bar{r}}{\bar{k}}.
\end{equation}
If this is the case, then the columns of the matrix
\begin{equation*}
	\sqrt{ \frac{ v }{ v-g } } ( I - E_s ) = \sqrt{ \frac{ v }{ v-g } } ( E_k + E_r )
\end{equation*}
also define a real ETF with $(M,N)=(v,v-g)$ in $V_k+V_r \, ( \cong \mathbb{R}^{v-g} )$, known as the \emph{Naimark complement} of $\{ \varphi_x : x \in X \}$.
We note that the spherical embedding of $\overline{\Gamma}$ consists of the normalized columns of $E_r$.

If $g=1$, then it follows from \eqref{angles} that $\Gamma$ is a complete bipartite graph, which we exclude here.
Hence we have $g\geqslant 2$.
Likewise, we have $f\geqslant 2$.
It follows that the real ETFs obtained in this way are always nontrivial, i.e., the pairs $(M,N)$ satisfy $2\leqslant N\leqslant M-2$.\footnote{The \emph{trivial} real ETFs are those with (a) $N=M$ (orthonormal bases); (b) $N=M-1$ (vertices of regular simplices); or (c) $N=1$. Cases (b) and (c) are the Naimark complements of each other.}

It is known (cf.~\cite[\S8.14]{BVM2022B}) that every dependent (i.e., $M>N$) real ETF gives rise to a \emph{regular two-graph}, that is to say, a $2$-design with block size three such that every $4$-set of vertices (or points) contains an even number of blocks.
Specifically, the vertices of the latter are indexed by the vectors of the former, where a $3$-set of vertices is a block if and only if an odd number of angles between the corresponding three vectors are obtuse angles.
Conversely, every regular two-graph is obtained in this way; see, e.g., \cite[\S\S10.3, 10.6]{BH2012B}, \cite[\S11.4]{GR2001B}.
By \cite[Proposition 10.3.2]{BH2012B}, one of the spherical embeddings of $\Gamma$ and $\overline{\Gamma}$ is a real ETF if and only if
\begin{equation}\label{regular two-graph}
	v=2(2k-\lambda-\mu)=2(2\bar{k}-\bar{\lambda}-\bar{\mu}).
\end{equation}

%%%%%%%%%%%%%%%%%%%%%
%%%%%%%%%%%%%%%%%%%%%
\section{Primitive rank \texorpdfstring{$3$}{3} graphs}
\label{sec: rank 3 graphs}

Brouwer and Van Maldeghem \cite{BVM2022B} described all \emph{primitive rank $3$ graphs}, i.e., those SRGs admitting primitive rank $3$ automorphism groups.
Recall that the \emph{socle} of a group $G$ is the subgroup generated by the minimal normal subgroups of $G$.

\begin{theorem}[{\cite[Theorem 11.1.1]{BVM2022B}}]
Let $\Gamma$ be a primitive strongly regular graph on $v$ vertices, and let $G$ be a primitive rank $3$ permutation group acting as an automorphism group of $\Gamma$.
Then one of the following holds.
\begin{enumerate}
\item[(i)] $T\times T \lhd G \leq T_0\wr 2$, where $T_0$ is a $2$-transitive group of degree $v_0$, the socle $T$ of $T_0$ is simple, and $v=v_0^2$.
\item[(ii)] The socle $L$ of $G$ is (nonabelian) simple.
\item[(iii)] The group $G$ is an affine group, i.e., $G$ has a regular elementary abelian normal subgroup and $v$ is a power of a prime.
\end{enumerate}
\end{theorem}

See also \cite{Fawcett2009M}.
Case (i) is discussed in \cite[\S11.2]{BVM2022B}.
There are $25$ kinds of groups, but all the groups generate the same class of the lattice graphs $L_2(q)$ (cf.~\cite[\S 1.1.8]{BVM2022B}).
Case (ii) has further four subcases: alternating socle, classical simple socle, exceptional simple socle, and sporadic simple socle; cf.~\cite[\S11.3]{BVM2022B}.
Case (iii) has further three subcases: infinite classes, extraspecial classes, and exceptional classes; cf.~\cite[\S11.4]{BVM2022B}.
The infinite families from Cases (ii) and (iii) are summarized in Tables \ref{Case (ii) families} and \ref{Case (iii) families}, respectively.\footnote{For some families, there are additional conditions to be of rank $3$.}
Most\footnote{The sporadic examples missing in \cite[\S11.5]{BVM2022B} are those in the classical simple socle subcase of Case (ii), i.e., Theorem 11.3.2\,(v), (vi), (vii), (x), and Theorem 11.3.3\,(iii) in \cite[\S11.3.2]{BVM2022B}. See also \url{https://homepages.cwi.nl/~aeb/math/srg/rk3/} .} of the sporadic examples are listed in \cite[\S11.5]{BVM2022B}.

\begin{table}[h]
\begin{minipage}{5.5cm}
\begin{center}
	\begin{tabular}{c|c}
	\hline\hline
	Graph& Reference \\
	\hline
	$T(n)$ & \cite[\S1.1.7]{BVM2022B} \\
	$\mathsf{Sp}_{2n}(q)$ & \cite[\S2.5]{BVM2022B} \\
	$\mathsf{O}_{2n+1}(q)$ & \cite[\S2.6]{BVM2022B} \\
	$\mathsf{O}_{2n}^+(q)$ & \cite[\S2.6]{BVM2022B} \\
	$\mathsf{O}_{2n+2}^-(q)$ & \cite[\S2.6]{BVM2022B} \\
	$\mathsf{U}_{2n}(\sqrt q)$ & \cite[\S2.7]{BVM2022B} \\
	$\mathsf{U}_{2n+1}(\sqrt q)$ & \cite[\S2.7]{BVM2022B} \\
	dual polar & \cite[\S2.2.11]{BVM2022B} \\
	half dual polar & \cite[\S2.2.12]{BVM2022B} \\
	$\mathit{NU}_n(q)$ & \cite[\S3.1.6]{BVM2022B} \\
	$\mathit{NO}_{2n}^{\varepsilon}(2)$ & \cite[\S3.1.2]{BVM2022B} \\
	$\mathit{NO}_{2n}^{\varepsilon}(3)$ & \cite[\S3.1.3]{BVM2022B} \\
	$\mathit{NO}_{2n+1}^{\varepsilon}(q)$ & \cite[\S3.1.4]{BVM2022B} \\
	$J_q(n,2)$ & \cite[\S3.5.1]{BVM2022B} \\
	$\mathsf E_{6,1}(q)$ & \cite[\S4.9]{BVM2022B} \\
	\hline\hline
\end{tabular}
\end{center}
\caption{Case (ii) families}\label{Case (ii) families}
\end{minipage}\!\!
\begin{minipage}{6cm}
\begin{center}
\begin{tabular}{c|c}
	\hline\hline
	Graph & Reference \\
	\hlineå
	$P(q)$ & \cite[\S1.1.9]{BVM2022B} \\
	$P^*(q)$ & \cite[\S7.3.6]{BVM2022B} \\
	Van Lint-Schrijver & \cite[\S7.3.1]{BVM2022B} \\
	$L_2(q)$ & \cite[\S1.1.8]{BVM2022B} \\
	$H_q(2,e)$ & \cite[\S3.4.1]{BVM2022B} \\
	$\mathit{VO}_{2n}^{\varepsilon}(q)$ & \cite[\S3.3.1]{BVM2022B} \\
	alternating forms & \cite[\S3.4.2]{BVM2022B} \\
	$\mathit{VD}_{5,5}(q)$ & \cite[\S3.3.3]{BVM2022B} \\
	$\mathit{VSz}(q)$ & \cite[\S3.3.1]{BVM2022B} \\
	\hline\hline
\end{tabular}
\end{center}
\caption{Case (iii) families}\label{Case (iii) families}
\end{minipage}
\end{table}

After many pages of cumbersome, case-by-case, and largely computer-assisted verifications of \eqref{equiangular} or \eqref{regular two-graph},\footnote{\label{details} Details are available at \url{https://zenodo.org/records/10849838}. We will use (*  *) to separate the cases. For instance, (*Lattice graph*)  denotes new cases of calculation of lattice graphs. } we obtained the following result:

\begin{theorem}\label{classification}
The primitive rank $3$ graphs whose spherical embeddings give rise to nontrivial real equiangular tight frames are given in Table \ref{classification table}.
\end{theorem}

\begin{table}[h]
\begin{center}
\begin{tabular}{c|c|c|c}
	\hline\hline
	Graph & $M$ & $N$ & $M-N$ \\
	\hline
	$\rule{0pt}{13pt} \mathit{NO}_{2n}^+(2)$\, $(n\geqslant 3)$ & $2^{n-1}(2^n-1)$ & $\frac{ ( 2^{n - 1 } - 1 ) ( 2^n - 1 ) }{ 3 }$ & $\frac{ 2^{2n}-1 }{ 3 }$ \\[1mm]
	$\overline{\mathit{NO}_{2n}^-(2)}$\, $(n\geqslant 2)$ & $2^{n-1}(2^n+1)$ & $\frac{ ( 2^{n - 1 } + 1 ) ( 2^n+ 1 ) }{ 3 }$ & $\frac{ 2^{2n}-1 }{ 3 }$ \\[1mm]
	$\mathit{NO}_{2n+1}^+(4)$\, $(n\geqslant 1)$ & $\frac{4^n ( 4^n + 1 )}{2}$ & $\frac{ 4^{2n}-1 }{ 3 }$ & $\frac{ (4^n+1)(4^n+2) }{ 6 }$ \\[1mm]
	$\overline{\mathit{NO}_{2n+1}^-(4)}$ \ $(n\geqslant 2)$ & $\frac{ 4^n ( 4^n - 1 ) }{2}$ & $\frac{ 4^{2n}-1 }{ 3 }$ & $\frac{ (4^n-1)(4^n-2) }{ 6 }$ \\[1mm]
	$\mathit{VO}_{2n}^+(2)$\, $(n\geqslant 2)$ & $2^{ 2n }$ & $2^{n- 1 } ( 2^n- 1 )$ & $2^{n-1}(2^n+1)$ \\[1mm]
	$\overline{\mathit{VO}_{2n}^-(2)}$\, $(n\geqslant 2)$ & $2^{ 2n }$ & $2^{n- 1 } ( 2^n+ 1 )$ & $2^{n-1}(2^n-1)$ \\[1mm]
	$\rule{0pt}{13pt} \overline{\mathsf{G}_2(2)} $ & $36$ & $21$ & $15$ \\
	$\rule{0pt}{13pt} \overline{\mathsf{M}_{22}} $ & $176$ & $154$ & $22$ \\
	\hline\hline
\end{tabular}
\end{center}
\caption{Primitive rank $3$ graphs yielding nontrivial real ETFs}\label{classification table}
\end{table}

\noindent
We also note the following isomorphisms: $T(5)\cong \overline{\mathit{NO}_4^-(2)}$ (cf.~\cite[\S10.3]{BVM2022B}), $\overline{T(8)}\cong \overline{\mathit{NO}_3^+(7)}\cong \mathit{NO}_6^+(2)$ (cf.~\cite[\S\S3.1.4, 3.6.1]{BVM2022B}), $\overline{L_2(4)}\cong \overline{\mathsf{O}_{4}^+(3)}\cong H_2(2,2)\cong \overline{\mathit{VO}_{2}^+(4)}\cong \mathit{VO}_4^+(2)$ (cf.~\cite[\S\S2.6.4, 3.3.1, 3.4.1]{BVM2022B}), $\mathit{NO}_5^-(3)\cong \overline{\mathit{NO}_6^-(2)}$ (cf.~\cite[\S10.15]{BVM2022B}), and $\mathit{VSz}(2)\cong\mathit{VO}^-_4(2)$ (cf.~\cite[\S2.5.5]{BVM2022B}).

%%%%%%%%%%%%%%%%%%%%%
%%%%%%%%%%%%%%%%%%%%%
\section{Rank \texorpdfstring{$3$}{3} graphs yielding equiangular tight frames}
\label{sec: real ETFs from SRGs}

In this section, we describe the primitive rank $3$ graphs found in Theorem \ref{classification} for the convenience of the reader.

%%%%%%%%%%%%%%%%%%%%%
\begin{example}[$\mathit{NO}^+_{ 2n } ( 2 )$]\label{example 1}
Equip $\mathbb{F}_2^{2n}$ with a nondegenerate quadratic form of Witt index $n$ and let $X$ be the set of nonsingular points.
Let $\Gamma=\mathit{NO}^+_{ 2n } ( 2 )$ have vertex set $X$, two vertices being adjacent when they are orthogonal.
See \cite[\S3.1.2]{BVM2022B}.
For $n\geqslant 3$, the graph $\mathit{NO}^+_{2n}(2)$ is a primitive rank $3$ graph with parameters $(2^{n-1}(2^n-1), 2^{2n-2}-1, 2^{2n-3}-2, 2^{n-2}(2^{n-1}+1) )$
and eigenmatrices
\begin{equation*}
	P = \begin{pmatrix} 1 & 2^{2n-2}-1 & 2^{n-1}(2^{n-1}-1) \\ 1 & 2^{ n- 2 } - 1 & - 2^{ n- 2 } \\ 1 & - 2^{ n- 1 } - 1 & 2^{ n- 1 } \end{pmatrix}\!, \ \ Q = \begin{pmatrix} 1 & \frac{ 2^{2n}-4 }{ 3 } & \frac{ ( 2^{ n- 1 } - 1 ) ( 2^n - 1 ) }{ 3 } \\[1mm] 1 & \frac{ 2^n-4 }{ 3 } & - \frac{ 2^n - 1 }{ 3 } \\[1mm] 1 & - \frac{ 2^n+2 }{ 3 } & \frac{ 2^n - 1 }{ 3 } \end{pmatrix}\!.
\end{equation*}
\end{example}

%%%%%%%%%%%%%%%%%%%%%
\begin{example}[$\overline{\mathit{NO}^-_{ 2n } ( 2 )}$]
Equip $\mathbb{F}_2^{2n}$ with a nondegenerate quadratic form of Witt index $n-1$ and let $X$ be the set of nonsingular points.
Let $\Gamma=\overline{\mathit{NO}^-_{ 2n } ( 2 )}$ have vertex set $X$, two vertices being adjacent when they are nonorthogonal.
See \cite[\S3.1.2]{BVM2022B}.
For $n\geqslant 2$, the graph $\overline{\mathit{NO}^-_{ 2n } ( 2 )}$ is a primitive rank $3$ graph with parameters $( 2^{n-1}(2^n+1), 2^{n-1}(2^{n-1}+1), 2^{n-2}(2^{n-1}+1), 2^{n-1}(2^{n-2}+1) )$ and eigenmatrices
\begin{equation*}
	P = \begin{pmatrix} 1 & 2^{n-1}(2^{n-1}+1) & 2^{2n-2}-1 \\ 1 & 2^{ n - 2 } & - 2^{ n - 2 } - 1 \\ 1 & - 2^{ n - 1 } & 2^{ n - 1 } - 1 \end{pmatrix}\!, \ \ Q = \begin{pmatrix} 1 & \frac{ 2^{2n}-4 }{ 3 } & \frac{ ( 2^{ n - 1 } + 1 ) ( 2^n + 1 ) }{ 3 } \\[1mm] 1 & \frac{ 2^n-2 }{ 3 } & - \frac{ 2^n + 1 }{ 3 } \\[1mm] 1 & - \frac{ 2^n+4 }{ 3 } & \frac{ 2^n + 1 }{ 3 } \end{pmatrix}\!.
\end{equation*}
\end{example}

%%%%%%%%%%%%%%%%%%%%%
\begin{example}[$\mathit{NO}^+_{ 2n + 1 } ( 4 )$]
Equip $\mathbb{F}_4^{2n+1}$ with a nondegenerate quadratic form and let $X$ be the set of nonsingular hyperbolic hyperplanes.
Let $\Gamma=\mathit{NO}^+_{ 2n+1 } ( 4 )$ have vertex set $X$, two vertices being adjacent when the restriction of the quadratic form to their intersection is degenerate.
See \cite[\S3.1.4]{BVM2022B}.
For $n\geqslant 1$, the graph $\mathit{NO}^+_{ 2n+1 } ( 4 )$ is a primitive rank $3$ graph with parameters $( 4^n ( 4^n + 1 )/2, ( 4^{ n - 1 } + 1 ) ( 4^n - 1 ), ( 4^{n-1}+2 ) ( 4^n-2 )/2, 4^n ( 4^{ n - 1 } + 1 )/2 )$ and eigenmatrices
\begin{gather*}
	P = \begin{pmatrix} 1 & ( 4^{ n - 1 } + 1 ) ( 4^n - 1 ) & 4^{n-1}(4^n-1) \\ 1 & 2 \cdot 4^{ n - 1 } - 1 & - 2 \cdot 4^{ n - 1 } \\ 1 & - 4^{ n - 1 } - 1 & 4^{ n - 1 } \end{pmatrix}\!, \\
	Q = \begin{pmatrix} 1 & \frac{ 2 ( 4^{ n - 1 } + 1 ) ( 4^n - 1 ) }{ 3 } & \frac{ 4^{2n}-1 }{ 3 } \\[1mm] 1 & \frac{ 4^n-2 }{ 3 } & - \frac{ 4^n + 1 }{ 3 } \\[1mm] 1 & - \frac{ 4 ( 4^{ n - 1 } + 1 ) }{ 3 } & \frac{ 4^n + 1 }{ 3 } \end{pmatrix}\!.
\end{gather*}
\end{example}

%%%%%%%%%%%%%%%%%%%%%
\begin{example}[$\overline{\mathit{NO}^-_{ 2n + 1 } ( 4 )}$]
Equip $\mathbb{F}_4^{2n+1}$ with a nondegenerate quadratic form and let $X$ be the set of nonsingular  elliptic hyperplanes.
Let $\Gamma=\overline{\mathit{NO}^-_{ 2n+1 } ( 4 )}$ have vertex set $X$, two vertices being adjacent when the restriction of the quadratic form to their intersection is nondegenerate.
See \cite[\S3.1.4]{BVM2022B}.
For $n\geqslant 2$, the graph $\overline{\mathit{NO}^-_{ 2n + 1 } ( 4 )}$ is a primitive rank $3$ graph with parameters $( 4^n ( 4^n - 1 )/2, 4^{ n - 1 } ( 4^n + 1 ), 4^n( 4^{n-1}+1 )/2, 4^{n-1} ( 4^n+2 )/2 )$ and 
eigenmatrices
\begin{gather*}
	P = \begin{pmatrix} 1 & 4^{ n - 1 } ( 4^n + 1 ) & ( 4^{ n - 1 } - 1 ) ( 4^n + 1 ) \\ 1 & 2 \cdot 4^{ n - 1 } & - 2 \cdot 4^{ n - 1 } - 1 \\ 1 & - 4^{ n - 1 } & 4^{ n - 1 } - 1 \end{pmatrix}\!, \\
	Q = \begin{pmatrix} 1 & \frac{ 2 ( 4^{ n - 1 } - 1 ) ( 4^n + 1 ) }{ 3 } & \frac{ 4^{2n}-1 }{ 3 } \\[1mm] 1 & \frac{ 4 ( 4^{ n - 1 } - 1 ) }{ 3 } & - \frac{ 4^n - 1 }{ 3 } \\[1mm] 1 & - \frac{ 4^n+2 }{ 3 } & \frac{ 4^n - 1 }{ 3 } \end{pmatrix}\!.
\end{gather*}
\end{example}

%%%%%%%%%%%%%%%%%%%%%
\begin{example}[$\mathit{VO}^+_{ 2n } ( 2 )$]
Equip $\mathbb{F}_2^{2n}$ with a nondegenerate quadratic form of Witt index $n$ and let $X=\mathbb{F}_2^{2n}$.
Let $\Gamma=\mathit{VO}^+_{2n} ( 2 )$ have vertex set $X$, two vertices being adjacent when their difference is isotropic.
See \cite[\S3.3.1]{BVM2022B}.
For $n\geqslant 2$, the graph $\mathit{VO}^+_{2n} (2)$ is a primitive rank $3$ graph with parameters $( 2^{ 2 n }, ( 2^{ n - 1 } + 1 ) ( 2^n - 1 ), ( 2^{ n - 1 } + 2 ) ( 2^{ n - 1 } -1 ), 2^{ n - 1 } ( 2^{ n - 1 } + 1 ) )$ and 
eigenmatrices
\begin{equation*}
	P = Q = \begin{pmatrix} 1 & ( 2^{ n - 1 } + 1 ) ( 2^n - 1 ) & 2^{ n - 1 } ( 2^n - 1 ) \\ 1 & 2^{ n - 1 } - 1 & - 2^{ n - 1 } \\ 1 & - 2^{ n - 1 } - 1 & 2^{ n - 1 } \end{pmatrix}\!.
\end{equation*}
\end{example}

%%%%%%%%%%%%%%%%%%%%%
\begin{example}[$\overline{\mathit{VO}^-_{ 2n } ( 2 )}$]\label{example 6}
Equip $\mathbb{F}_2^{2n}$ with a nondegenerate quadratic form of Witt index $n-1$ and let $X=\mathbb{F}_2^{2n}$.
Let $\Gamma=\overline{\mathit{VO}^-_{2n} ( 2 )}$ have vertex set $X$, two vertices being adjacent when their difference is nonisotropic.
See \cite[\S3.3.1]{BVM2022B}.
For $n\geqslant 2$, the graph $\overline{\mathit{VO}^-_{ 2n } ( 2 )}$ is a primitive rank $3$ graph with parameters $( 2^{ 2 n }, 2^{ n - 1 } ( 2^n + 1 ), 2^{n-1}( 2^{ n - 1 } + 1 ), 2^{ n - 1 } ( 2^{ n - 1 } +1 ) )$ and 
eigenmatrices
\begin{gather*}
	P = \begin{pmatrix} 1 & 2^{ n - 1 } ( 2^n + 1 ) & ( 2^{ n - 1 } - 1 ) ( 2^n + 1 ) \\ 1 & 2^{ n - 1 } & - 2^{ n - 1 } - 1 \\ 1 & - 2^{ n - 1 } & 2^{ n - 1 } - 1 \end{pmatrix}\!, \\
	Q = \begin{pmatrix} 1 & ( 2^{ n - 1 } - 1 ) ( 2^n + 1 ) & 2^{ n - 1 } ( 2^n + 1 ) \\ 1 & 2^{ n - 1 } - 1 & - 2^{ n - 1 } \\ 1 & - 2^{ n - 1 } - 1 & 2^{ n - 1 } \end{pmatrix}\!.
\end{gather*}
\end{example}

%%%%%%%%%%%%%%%%%%%%%
\begin{example}[$\overline{\mathsf{G}_2(2)}$]
Let $X$ be the set consisting of the $7$ points, $7$ lines, and $21$ flags of the Fano plane, together with an additional element denoted by $\infty$.
The graph $\overline{\mathsf{G}_2(2)}$ has vertex set $X$ and the adjacency is defined as follows:
The vertex $\infty$ is adjacent to the flags.
The points form a clique, and so do the lines.
A point and a line are adjacent when they are incident.
A point (resp.~line) is adjacent to the flags whose lines (resp.~ points) are not on it.
Finally, two flags are adjacent when either they are not disjoint, or they are disjoint and the point of each of them is not on the line of the other.
See \cite[\S10.14]{BVM2022B}.
The graph $\overline{\mathsf{G}_2(2)}$ is a primitive rank $3$ graph with parameters $(36,21,12,12)$ and eigenmatrices
\begin{equation*}
	P=\begin{pmatrix} 1 & 21 & 14 \\ 1 & 3 & -4 \\ 1 & -3 & 2 \end{pmatrix}\!, \quad Q=\begin{pmatrix} 1 & 14 & 21 \\[1mm] 1 & 2 & -3 \\[1mm] 1 & -4 & 3 \end{pmatrix}\!.
\end{equation*}
\end{example}

%%%%%%%%%%%%%%%%%%%%%
\begin{example}[$\overline{\mathsf{M}_{22}}$]
Let $X$ be the set of $176$ blocks of the unique quasi-symmetric $2$-$(22,7,16)$ design with block intersection numbers $1$ and $3$.
Let $\Gamma=\overline{\mathsf{M}_{22}}$ have vertex set $X$, two vertices being adjacent when they intersect in $3$ points.
See \cite[\S10.51]{BVM2022B}.
The graph $\overline{\mathsf{M}_{22}}$ is a primitive rank $3$ graph with parameters $(176,105,68,54)$ and eigenmatrices
\begin{equation*}
	P=\begin{pmatrix} 1 & 105 & 70 \\ 1 & 17 & -18 \\ 1 & -3 & 2 \end{pmatrix}\!, \quad Q=\begin{pmatrix} 1 & 21 & 154 \\[1mm] 1 & \frac{17}{5} & -\frac{22}{5} \\[1mm] 1 & -\frac{27}{5} & \frac{22}{5} \end{pmatrix}\!.
\end{equation*}
\end{example}

\begin{remark}
The six infinite families of SRGs in Examples \ref{example 1}--\ref{example 6} are part of more general families of association schemes obtained from actions of classical groups, and their eigenmatrices were extensively studied.
See, e.g., \cite{Bannai1990P,BHS1990JCTA,BKS1990MFSKUA,BSHW1991JA,Kwok1991GC,Kwok1992EJC,ST2006P,Tanaka2001M,Tanaka2002EJC,Tanaka2004AG}.
It may be interesting to point out that these SRGs and association schemes were used to construct Ramanujan graphs; cf.~\cite{BST2004EJC,BST2009DM}.
\end{remark}

Observe that some of the real ETFs in Theorem \ref{classification} have the same parameters up to Naimark complements (see also Table \ref{comparison} in Section \ref{sec: Waldron ETFs}).
We wonder if, for example, the real ETFs obtained from $\mathit{NO}_{2n+1}^+(4)$ and $\overline{\mathit{NO}_{4n}^-(2)}$ are equivalent, i.e., the graphs $\mathit{NO}_{2n+1}^+(4)$ and $\overline{\mathit{NO}_{4n}^-(2)}$ are switching equivalent.

We also note that the real ETFs in Theorem \ref{classification} are still realized only as Gram matrices since the isomorphism $V_s\cong\mathbb{R}^N$ (where $N=g$) is not canonical.
However, the families $\mathit{VO}^+_{ 2n } ( 2 )$ and $\overline{\mathit{VO}^-_{ 2n } ( 2 )}$ are Cayley graphs on elementary abelian $2$-groups, and their real ETFs can be given as explicit $N$-dimensional real unit column vectors.
For $\mathit{VO}^+_{ 2n} ( 2 )$, we may assume that the quadratic form is given by
\begin{equation*}
	q(x)=x_1x_2+\cdots+x_{2n-1}x_{2n}
\end{equation*}
for $x=(x_1,\dots,x_{2n})\in\mathbb{F}_2^{2n}$.
The associated bilinear form is then given by
\begin{equation*}
	B(x,y)=(x_1y_2+x_2y_1)+\cdots+(x_{2n-1}y_{2n}+x_{2n}y_{2n-1})
\end{equation*}
for $x=(x_1,\dots,x_{2n}),y=(y_1,\dots,y_{2n})\in\mathbb{F}_2^{2n}$.
Let $X_1$ be the set of nonisotropic points in $X=\mathbb{F}_2^{2n}$.
Then $|X_1|=g=N=2^{n-1}(2^n-1)$, and the following column vectors form an orthonormal basis of $V_s$:
\begin{equation*}
	\frac{1}{2^n}\big((-1)^{B(x,z)}\big)_{x\in X} \quad (z\in X_1).
\end{equation*}
See, e.g., \cite{Tanaka2004AG}.
Hence, with respect to this basis, the vectors $\varphi_x$ $(x\in X)$ in the real ETF are expressed as $2^{n-1}(2^n-1)$-dimensional column vectors
\begin{equation*}
	\frac{1}{\sqrt{2^{n-1}(2^n-1)}}\big((-1)^{B(x,z)}\big)_{z\in X_1} \quad (x\in X).
\end{equation*}
A similar discussion applies to $\overline{\mathit{VO}^-_{ 2n } ( 2 )}$.

%%%%%%%%%%%%%%%%%%%%%
%%%%%%%%%%%%%%%%%%%%%
\section{Real equiangular tight frames having rank \texorpdfstring{$3$}{3} descendants}
\label{sec: Waldron ETFs}

We mentioned earlier that every dependent real ETF corresponds to a regular two-graph.
If, moreover, the former is nontrivial, then the derived designs of the latter are SRGs with $v=M-1$ and $k=2\mu$.
The derived designs are called the \emph{descendants} \cite[\S10.2]{BH2012B} and also the \emph{neighborhoods} \cite[\S11.5]{GR2001B} of the latter.
Here, let us temporarily call them the \emph{descendants of the real ETF}.
Conversely, if $\Gamma$ is an SRG with $k=2\mu$, then the matrix
\begin{equation*}
	I +\frac{1}{1+2r} \left(\begin{array}{c|c} 0 & \bm{1}^{\mathsf{T}} \\ \hline \bm{1} & J-I-2A \end{array}\right)
\end{equation*}
is the Gram matrix of a nontrivial real ETF with $M=v+1$ and $N=g+1$, where $\bm{1}$ denotes the all one's vector.
The graphs $\Gamma$ and $\overline{\Gamma}$ are then descendants of this real ETF and its Naimark complement, respectively.
See \cite[\S10.3]{BH2012B}, \cite[\S11.6]{GR2001B}, and also \cite[\S5]{Waldron2009LAA}.
Our findings are as follows:\footnote{See Footnote \ref{details}.}

\begin{theorem}\label{rank 3 Waldron}
The primitive rank $3$ graphs, up to taking complements, that are descendants of nontrivial real equiangular tight frames, are given in Table \ref{Waldron table}.
\end{theorem}

\begin{table}[ht]
\begin{center}
\begin{tabular}{c|c|c|c}
	\hline\hline
	Graph & $M$ & $N$ & $M-N$ \\
	\hline
	$\rule{0pt}{13pt} \mathsf{Sp}_{2n}(2)$\, $(n\geqslant 2)$ & $2^{2n}$ & $2^{n-1}(2^n-1)$ & $2^{n-1}(2^n+1)$ \\[1mm]
	$\mathsf{O}_{2n}^+(2)$\, $(n\geqslant 2)$ & $2^{n-1}(2^n+1)$ & $\frac{ 2^{2n}-1 }{ 3 }$ & $\frac{ ( 2^{ n - 1 } + 1 ) ( 2^n + 1 ) }{ 3 }$ \\[1mm]
	$\mathsf{O}_{2n+2}^-(2)$\, $(n\geqslant 2)$ & $2^n(2^{n+1}-1)$ & $\frac{ (2^n-1)(2^{n+1}-1) }{ 3 }$ & $\frac{ 2^{2n+2}-1 }{ 3 }$ \\[1mm]
	$P(q)$ & $q+1$ & $\frac{ q+1 }{2}$ & $\frac{ q+1 }{2}$ \\[1mm]
	$P^*(q)$ & $q+1$ & $\frac{ q+1 }{2}$ & $\frac{ q+1 }{2}$ \\[1mm]
	McLaughlin & $276$ & $23$ & $253$ \\[1mm]
	$P^{**}(529)$ & $530$ & $265$ & $265$ \\[1mm]
	$(2209,1104,551,552)$ & $2210$ & $1105$ & $1105$ \\
	\hline\hline
\end{tabular}
\end{center}
\caption{Primitive rank $3$ descendants of nontrivial real ETFs\\[2mm] \small For $P(q)$, $q$ is a prime power and $q\equiv 1\, (\operatorname{mod}\,4)$. For $P^*(q)$, $q$ is an even power of a prime $p\equiv 3\, (\operatorname{mod}\,4)$. For the McLaughlin graph, see \cite[\S10.61]{BVM2022B}. The graph $P^{**}(529)$ is the sporadic Peisert graph \cite[\S10.70]{BVM2022B}. There are precisely three rank $3$ graphs with the last parameters, two of which are $P(2209)$ and $P^*(2209)$; cf.~\cite[\S10.86]{BVM2022B}.}\label{Waldron table}
\end{table}

\noindent
We also note the following isomorphisms: $\mathsf{O}_{2n+1}(2)\cong\mathsf{Sp}_{2n}(2)$ (cf.~\cite[\S2.6.3]{BVM2022B}), $\overline{T(6)}\cong\mathit{NO}^-_4(3)\cong\overline{\mathit{NO}^+_3(5)}\cong\mathsf{Sp}_4(2)$ (cf.~\cite[\S\S3.1.3, 3.1.4, 10.5]{BVM2022B}), $J_2(4,2)\cong\mathsf{O}^+_6(2)$ (cf.~\cite[\S10.13]{BVM2022B}), and $L_2(3)\cong \mathit{VO}_2^+(3)\cong P(9)\cong\mathsf{O}_4^+(2)$ (cf.~\cite[\S10.2]{BVM2022B}).
It should be mentioned that the graph $\mathit{NU}_3(3)$ satisfies $k=2\mu$ but is of rank $4$; cf.~\cite[\S10.22]{BVM2022B}.
We do not describe the graphs in Table \ref{Waldron table} to keep the paper concise.
See the references given.

The Paley graphs $P(q)$ and the Peisert graphs $P^*(q)$ give rise to real ETFs having the same parameters, so we ask if these real ETFs are equivalent, i.e., the disjoint union $K_1+P(q)$ is switching equivalent to $K_1+P^*(q)$.
We may ask the same question for $P^{**}(529)$ and the last graph in Table \ref{Waldron table}.

\begin{table}[ht]
\begin{center}
\begin{tabular}{c|c|c|c}
	\hline\hline
	Thm.~\ref{classification} & Thm.~\ref{rank 3 Waldron} & $M$ & $\{N,M-N\}$ \\
	\hline
	$\rule{0pt}{13pt} \mathit{NO}_{2n}^+(2)$ & $\mathsf{O}_{2n}^-(2)$ & $2^{n-1}(2^n-1)$ & $\left\{\frac{ ( 2^{ n - 1 } - 1 ) ( 2^n - 1 ) }{ 3 },\frac{ 2^{2n}-1 }{ 3 }\right\}$ \\[2mm]
	$\overline{\mathit{NO}_{2n}^-(2)}$ & $\mathsf{O}_{2n}^+(2)$ & $2^{n-1}(2^n+1)$ & $\left\{\frac{ ( 2^{ n - 1 } + 1 ) ( 2^n + 1 ) }{ 3 },\frac{ 2^{2n}-1 }{ 3 }\right\}$ \\[2mm]
	$\mathit{NO}_{2n+1}^+(4)$ & $\mathsf{O}_{4n}^+(2)$ & $2^{2 n-1} \left(2^{2n}+1\right) $ & $\left\{\frac{2^{4n}-1}{3} , \frac{(2^{2n-1}+1)(2^{2n}+1)}{3} \right\}$ \\[2mm]	
	$\overline{\mathit{NO}_{2n+1}^-(4)}$ & $\mathsf{O}_{4n}^-(2)$ & $2^{2 n-1} \left(2^{2n}-1\right) $ & $\left\{\frac{2^{4n}-1}{3}  ,\frac{(2^{2n-1}-1)(2^{2n}-1)}{3}  \right\}$ \\[2mm]	
	$\mathit{VO}_{2n}^+(2),\overline{\mathit{VO}_{2n}^-(2)}$ & $\mathsf{Sp}_{2n}(2)$ & $2^{ 2 n }$ & $\big\{2^{ n - 1 } ( 2^n - 1 ),2^{n-1}(2^n+1)\big\}$ \\[1mm]
	$\overline{\mathsf{G}_2(2)}$ & $\mathsf{O}_6^+(2)$ & $36$ & $\{21,15\}$ \\[1mm]
	\hline\hline
\end{tabular}
\end{center}
\caption{Comparison of the parameters of real ETFs}\label{comparison}
\end{table}

Some of the real ETFs found in Theorems \ref{classification} and \ref{rank 3 Waldron} also have the same parameters; see Table \ref{comparison}.
We again wonder if some of these real ETFs are equivalent.
If the spherical embedding of an SRG $\Gamma$ is a real ETF, then its descendants are obtained by first removing a vertex from $\Gamma$ and then switching it with respect to the neighbors and the nonneighbors of the removed vertex.
For example, if $\Gamma=\mathit{NO}_{2n}^+(2)$, then we ask if the resulting graph is isomorphic to $\mathsf{O}_{2n}^-(2)$.
It seems that this problem was already addressed in \cite[Remark 4.6]{FIJK2021DCC} for the first, second, and fifth rows in Table \ref{comparison}.
We may remark that if $\Gamma$ has parameters $(v,k,\lambda,\mu)$ (cf.~\eqref{regular two-graph}), then the parameters of the descendants are given by $(v-1,2(k-\mu),k+\lambda-2\mu,k-\mu)$; see \cite[Proposition 10.3.2]{BH2012B}.

%%%%%%%%%%%%%%%%%%%%%
%%%%%%%%%%%%%%%%%%%%%
\section*{Acknowledgments}

Eiichi Bannai and Etsuko Bannai thank National Center for Theoretical Sciences (NCTS) at National Taiwan University (NTU) for inviting them to stay for 14 months (Dec.~2020--Jan.~2022).
The last stage of this work was done while they were visiting there.
Hajime Tanaka was supported by JSPS KAKENHI Grant Number JP20K03551.
Wei-Hsuan Yu and Chin-Yen Lee were supported by MOST under Grant MOST109-2628-M-008-002-MY4.

\end{document}